\documentclass[p4paper, 12pt]{article}
\usepackage{epsfig, color, epstopdf}
\usepackage{amssymb,amsmath,amsthm,latexsym,lineno}
\usepackage{graphicx,graphics,amsfonts,indentfirst}
\usepackage{mathrsfs,eufrak, bbm, dsfont,yfonts}
\usepackage{url,tikz,yhmath,appendix,booktabs}
\usepackage{multirow}
\usepackage{kantlipsum}
\allowdisplaybreaks
\usepackage{setspace}

\usepackage{stmaryrd}
\newcommand \brk[1]
{
	\llbracket #1\rrbracket
}

\newtheorem{theo}{Theorem}
\newtheorem{lem}[theo]{Lemma}
\newtheorem{pro}[theo]{Proposition}
\newtheorem{prob}{Problem}
\newtheorem{exa}[theo]{Example}
\newtheorem{con}{Conjecture}
\newtheorem{cor}[theo]{Corollary}
\newtheorem{clm}{Claim}

\usetikzlibrary{backgrounds}
\usetikzlibrary{arrows}
\usetikzlibrary{shapes,shapes.geometric,shapes.misc}
\pgfdeclarelayer{edgelayer}
\pgfdeclarelayer{nodelayer}
\pgfsetlayers{background,edgelayer,nodelayer,main}
\tikzstyle{none}=[inner sep=0mm]
\tikzstyle{every loop}=[]

\tikzstyle{dotted}=[dash pattern=on \pgflinewidth off 2pt]
\tikzstyle{dashed}=[dash pattern=on 3pt off 3pt]

\newcommand \tikzp[2]
{
\begin{center}
\begin{tikzpicture}[scale=#1]
#2
\end{tikzpicture}
\end{center}
}

\tikzstyle{new style 0}=[fill=black, draw=black, shape=circle]
\tikzstyle{red style 1}=[fill=red, draw=black, shape=circle]
\tikzstyle{blue style 2}=[fill=blue, draw=black, shape=circle]
\tikzstyle{white style 4}=[fill=white, draw=black, shape=circle]
\tikzstyle{bklack style 5}=[fill=black, draw=black, shape=rectangle]
\tikzstyle{red style 3}=[fill=red, draw=black, shape=rectangle]
\tikzstyle{yellow style 7}=[fill=yellow, draw=black, shape=rectangle]
\tikzstyle{new style 8}=[fill={rgb,255: red,0; green,132; blue,0}, draw={rgb,255: red,0; green,131; blue,0}, shape=circle]

\tikzstyle{new edge style 0}=[-]
\tikzstyle{new edge style 1}=[-, draw=red]
\tikzstyle{new edge style 2}=[-, draw=blue]
\tikzstyle{new edge style 3}=[-, draw={rgb,255: red,0; green,156; blue,0}]

\tikzstyle{cblue}=[circle, draw, thin,fill=blue!20, scale=0.5]

\newcommand \equ[2]
{
\begin{equation}\label{#1}
#2
\end{equation}
}

\newcommand \eqn[2]
{
\begin{eqnarray}\label{#1}
#2
\end{eqnarray}
}

\newcommand \prop[2]
{\begin{pro}
\label{#1} #2
\end{pro}
}

\newcommand \prom[2]
{\begin{prob}
\label{#1} #2
\end{prob}
}

\newcommand \corr[2]
{\begin{cor}
\label{#1} #2
\end{cor}
}

\newcounter{countcase}

\newcounter{countclaim}

\def \proof {\noindent {{\it Proof}}.\setcounter{countcase}{0} \setcounter{clm}{0}}

\newcommand{\proofn}{{\noindent {\em Proof}.}\setcounter{countcase}{0} }

\newcommand{\claimend}{{\hfill $\natural $}}

\newcommand{\proofend}{{\hfill$\Box$}}

\def \N {{\mathbb N}}

\def \setc {{\cal C}}

\def \seti {{\cal I}}
\def \seth {{\cal H}}
\def \cone {{\cal E}_{e}}

\def \scre {{\mathscr E}}
\def \scrs {{\mathscr S}}
\def \scrh {{\mathscr H}}
\def \scrx {{\mathscr X}}

\def \gt {{\cal G}}
\def \gtp {{\cal GT}}

\def \DPg {{\cal DP}^*}

\newcommand \dpg[1]
{
P_{DP}(#1)\approx P(#1)
}

\newcommand \dpl[1]
{
P_{DP}(#1)<  P(#1)
}

\newcommand \spnn[2]
{
#1\langle #2 \rangle
}

\newcommand \spng[2] 
{
H[#1]|_{#2}
}

\begin{document}

\baselineskip 0.6 cm

\title{DP color functions versus  chromatic polynomials}

\author{
Fengming Dong\thanks{Corresponding author. 
Email: fengming.dong@nie.edu.sg
and donggraph@163.com.}
\\
\small National Institute of Education,
Nanyang Technological University, Singapore
\\
Yan Yang\thanks{Email:  yanyang@tju.edu.cn}
\\
\small School of Mathematics, Tianjin University, China
}
\date{}

\maketitle{}

\begin{abstract}
For any graph $G$, the chromatic polynomial of $G$ is the function
$P(G,m)$ which counts the number of 
proper $m$-colorings of $G$
for each positive integer $m$.
The DP color function $P_{DP}(G,m)$ of $G$, 
introduced by 
Kaul and Mudrock in 2019, \label{oc1}
is a generalization of $P(G,m)$ with 
$P_{DP}(G,m)\le P(G,m)$ for each positive integer $m$.
Let $P_{DP}(G)\approx P(G)$ 
(resp. $P_{DP}(G)<  P(G)$)
denote the property that 
$P_{DP}(G,m)=P(G,m)$ 
(resp. $P_{DP}(G,m)<P(G,m)$)
holds for sufficiently large integers $m$.
It is an interesting problem of finding 
graphs $G$ for which  $P_{DP}(G)\approx P(G)$ 
(resp. $P_{DP}(G,m)<P(G,m)$) holds.
Kaul and Mudrock showed that 
if $G$ has an even girth, then 
$P_{DP}(G)<P(G)$
and Mudrock and Thomason recently proved that 
$P_{DP}(G)\approx P(G)$ holds for each graph $G$ 
which has a dominating vertex. 
We shall generalize their results in this article.
For each edge $e$ in $G$, let $\ell(e)=\infty$
if $e$ is a bridge of $G$, and 
let $\ell(e)$ be the length of a shortest cycle in $G$ 
containing $e$ otherwise.
We first show that if 
$\ell(e)$ is even for some edge $e$ in $G$,
then $P_{DP}(G)<P(G)$ holds.
However, the converse statement of this conclusion fails
with infinitely many counterexamples.  \label{cr1}
We then prove that $P_{DP}(G)\approx P(G)$ holds 
for every graph $G$ that  \label{cr2}
contains a spanning tree $T$ such that 
for each $e\in E(G)\setminus E(T)$, 
$\ell(e)$ is odd and $e$ is contained in a cycle $C$ 
of length $\ell (e)$ with the property that  
$\ell(e')<\ell(e)$ for each 
$e'\in E(C)\setminus (E(T)\cup \{e\})$.
Some open problems are proposed in this article.
\end{abstract} 




\section{Introduction}
\label{sec1}

In this article, we consider simple graphs only,
unless otherwise stated. 
For any graph $G$, let $V(G)$ and $E(G)$ be 
its vertex set and edge set respectively.
For any nonempty subset $S$ of $V(G)$,
let $G[S]$ be the subgraph of $G$ induced by $S$, 
i.e., the subgraph with vertex set $S$ and 
edge set $\{uv\in E(G):u,v\in S\}$,
where $uv$ denotes the edge 
joining $u$ and $v$,  \label{cr3}
and let $G-S$ be the subgraph $G[V(G)\setminus S]$ when $S\ne V(G)$.
In particular, if $S=\{v\}$ for $v\in V(G)$, 
write $G-v$ for $G-S$.
For $A\subseteq E(G)$, 
let $\spnn{G}{A}$ denote the spanning subgraph
of $G$ with edge set $A$,
and let $G-A$ be $\spnn{G}{E(G)\setminus A}$.
In particular, for $e\in E(G)$, 
$G-\{e\}$ is written as $G-e$.
For two disjoint subsets $S_1$ and $S_2$ of $V(G)$,
let $E_G(S_1,S_2)$ (or simply $E(S_1,S_2)$)  denote the set 
$\{uv\in E(G): u\in S_1,v\in S_2\}$.
For any $u\in V(G)$, let 
$N_G(u)$ (or simply $N(u)$) be the set of neighbors of $u$ in $G$
and $d_G(u)$ (or simply $d(u)$) be the degree of $u$ 
in $G$.
The reader may refer to \cite{BM76}
for other terminology and notation.

\subsection{Proper coloring, list coloring
and DP coloring}

Let $\N$ denote the set of positive integers.
For any $n\in \N$, let $\brk{n}=\{1,2,\cdots,n\}$.
For any graph $G$ and  $k\in \N$,
a {\it proper $k$-coloring} of $G$ is a mapping 
of $f:V(G)\rightarrow \brk{k}$
such that $f(u)\ne f(v)$ for each edge $uv\in E(G)$.
The {\it chromatic polynomial} $P(G,k)$ of $G$, 
introduced by Birkhoff~\cite{birk} in 1912, 
is the function which counts the number of 
proper $k$-colorings of $G$ for each $k\in \N$.
Note that $P(G,k)$ is indeed a polynomial in $k$ for each 
$k\in \N$ (see \cite{birk2, dong0, rea1, Whi1932}).
The {\it chromatic number} of $G$, denoted by $\chi(G)$,
is the minimum number $k\in \N$ such that 
$G$ admits a proper $k$-coloring. 
Obviously, $\chi(G)$ is the minimum number $k\in \N$ 
such that $P(G,k)>0$.
For more details on chromatic polynomials,
we refer the readers to~\cite{birk,
birk2, dong0, dong1, Jackson2015, rea1, rea2,
Royle2009, Whi1932}.

List coloring was introduced by Vizing \cite{viz} and 
Erd\H{o}s, Rubin and Taylor \cite{erdos} independently. 
A list coloring of $G$ is associated with 
a {\it list assignment} $L$, where  
$L(v)$ is a subset of $\N$ for each $v\in V(G)$.
Given a list assignment $L$ of $G$,
a {\it proper} $L$-coloring of $G$ 
is a mapping $f:V(G)\rightarrow \N$ 
such that $f(v)\in L(v)$ for each $v\in V(G)$
and $f(u)\ne f(v)$ for each edge $uv\in E(G)$.
If $L(v)=\brk{k}$ for each $v\in V(G)$,
then a proper $L$-coloring of $G$ is 
a proper $k$-coloring of $G$. 
If $|L(v)| = m$ for each $v \in V (G)$, then $L$ is called an {\it $m$-assignment} .
The {\it list chromatic number} of $G$,
denoted by $\chi_l(G)$, is the smallest $m$ 
such that $G$ has a proper $L$-coloring 
for every $m$-assignment $L$ of $G$. 
By definition, $\chi(G)\le \chi_l(G)$.
Due to  Noel, Reed and Wu~\cite{Noel}, 
$\chi(G)=\chi_l(G)$ holds whenever $\chi(G)\ge (|V(G)|-1)/2$.

For each list-assignment $L$ of $G$,
let $P(G,L)$ be the number of proper $L$-colorings.
For each $m\in \N$, 
let $P_l(G,m)$ be the minimum value of $P(G,L)$ 
among all $m$-assignments $L$.
We call $P_l(G,m)$ the {\it list color function} of $G$.
By definition,  $P_l(G,m)\le P(G,m)$ for each $m\in \N$.
Wang, Qian and Yan~\cite{wang} showed that $P_l(G,m)= P(G,m)$ holds when 
$G$ is connected and $m>(|E(G)|-1)/\ln{(1+\sqrt 2)}$. 
The survey by Thomassen~\cite{Thomassen1} 
provided some known results
and open questions on the list color function.

DP-coloring was introduced by Dvor\'ak and Postle
\cite{dp} for the purpose of proving that 
every planar graph without cycles of lengths 4 to 8 is
3-choosable. 
DP-coloring is a generalization of list coloring,
and a formal definition is given below.
For a graph $G$, a {\it cover} of $G$ is an ordered pair 
$\seth = (L, H)$, 
where $H$ is a graph and 
$L$ is a mapping from $V(G)$ to the power set of $V(H)$
satisfying the four conditions below:
\begin{enumerate}
\item 
the sets $\{L(u):u\in V (G)\}$ 
is a partition of 
$V(H)$ of size $|V(G)|$; \label{oc2}

\item for every $u\in V (G)$, 
$H[L(u)]$ is a complete graph;

\item if $u$ and $v$ are non-adjacent vertices in $G$,
then $E_H (L(u), L(v))=\emptyset$; and 

\item for each edge $uv\in E(G)$,  
$E_H (L(u), L(v))$ is a matching.
\end{enumerate} \label{cr8}

An {\it $\seth$-coloring} of $G$ is 
an independent set $I$ of $H$ with $|I|=|V(G)|$. 
Clearly, for each $\seth$-coloring $I$ of $G$, 
$|I\cap L(u)| = 1$ holds for each $u \in V (G)$. 
A cover $\seth = (L, H)$  of $G$ is called 
an {\it $m$-fold} cover if 
$|L(u)| = m$ for each $u\in V (G)$.
The {\it $DP$-chromatic number} of $G$,
denoted by $\chi_{DP}(G)$, 
is the minimum integer $m$ such that 
$G$ has a $\seth$-coloring for every $m$-fold cover 
$\seth=(L,H)$.
By definition, $\chi(G)\le \chi_l(G)\le \chi_{DP}(G)$.
Bernshteyn, Kostochka and Zhu~\cite{Bern} showed that  
for any $n\in \N$,
if $r(n)$ is the minimum number $r\in \N$ 
such that $\chi(G)
=\chi_{DP}(G)$ holds 
for every $n$-vertex graph $G$ with $\chi(G)\ge r$,
then $n-r(n)=\Theta(\sqrt n)$. 

For any cover $\seth$ of $G$, let $P_{DP}(G,\seth)$ be the number of
$\seth$-colorings of $G$.
For each $m\in \N$,
let $P_{DP}(G,m)$ be the minimum value 
of $P_{DP}(G,\seth)$ among all $m$-fold covers $\seth$ of $G$.
We call $P_{DP}(G,m)$ the {\it DP color function} of $G$,
which was introduced by Kaul and Mudrock~\cite{kaul1}.
For any $m$-assignment $L$ of $G$, 
$P(G,L)=P_{DP}(G,\seth)$ holds for the $m$-fold cover $\seth=(L',H)$,
where $L'(v)=\{(v,j): j\in L(v)\}$ for each $v\in V(G)$
and for each edge $uv\in E(G)$,
$E_H(L'(u),L'(v))=\{(u,j)(v,j): j\in \N, (u,j)\in L'(u),
(v,j)\in L'(v)\}$. 
Thus, $P_{DP}(G,m)\le P_{l}(G,m)\le P(G,m)$ holds for each $m\in \N$. \label{cr7}

\subsection{Main results
\label{sec1-2}}

For any graph $G$, by definition,
$P_{DP} (G, m) \le P(G, m)$ 
holds for all integers $m\in N$.
Thus, exactly one of the 
following three  properties
holds:
\begin{enumerate}
\item there exists $N\in \N$ such that 
$P_{DP} (G, m) = P(G, m)$ for all integers $m \ge N$;
\item there exists $N\in \N$ such that 
$P_{DP} (G, m) < P(G, m)$ for all integers $m \ge N$;
and
\item there exist two infinite sets
$\{m_i\in \N:i\in \N\}$ and $\{n_i\in \N:i\in \N\}$
such that for all $i\in \N$, 
both $P_{DP} (G, m_i)  =P(G, m_i)$  
and $P_{DP} (G, n_i)  <P(G, n_i)$ hold.
\end{enumerate}

Two questions proposed by Kaul and Mudrock \cite{kaul1}  
are closed related to property (iii),
and there would be no graphs having property (iii)
if the answer to any one of them had been yes.
Question 7 in \cite{kaul1} asks 
if, for any graph $G$, there always exist an 
$N\in \N$ and a polynomial $p(m)$ 
such that $P_{DP}(G, m) =p(m)$ whenever $m \ge N$.
Halberg, Kaul, Liu, Mudrock,
Shin and Thomason~\cite{Halb} showed that 
this question has a positive answer for 
each graph $G$ with a vertex $v$ such that 
$G-v$ is acyclic.
Question 15 in \cite{kaul1} asks 
if $P_{DP}(G, m_0) =P(G, m_0)$ for some 
$m_0\ge \chi(G)$ implies that  
$P_{DP}(G, m) =P(G, m)$ for all $m \ge m_0$.
Unfortunately,  \label{oc3}
Bui, Kaul, Maxfield, Mudrock, Shin and Thomason~\cite{Bui}
found graphs with negative answer to the second question.

For any one of the above properties,
it is an interesting problem 
of knowing which graphs have this property.
For convenience purposes, 
let $\dpg{G}$ (resp. $\dpl{G}$)
denote property (i) (resp. property (ii)) above
for a graph $G$.

\prom{prom0-1}
{
Is it true that for each graph $G$, either $\dpg{G}$ or $\dpl{G}$?
}

So far the comparison of DP color functions 
with chromatic polynomials 
focuses on following problem.

\prom{prom0}
{
Determine the set of graphs $G$ such that 
$\dpg{G}$ holds and the set of graphs $G$ such that
$\dpl{G}$ holds.  
}

Kaul and Mudrock~\cite{kaul1} obtained some important  results
on the study of Problem~\ref{prom0}.
For example, they showed that if there exists an edge 
$e$ in $G$ such that $P(G-e,m)<mP(G,m)/(m-1)$, then 
$P_{DP}(G,m)<P(G,m)$ holds (see Theorem~\ref{th3-1}).

For each edge $e$ in $G$, if $e$ is a bridge of $G$,
let $\ell_G(e)=\infty$;
otherwise, let $\ell_G(e)$ be 
the length of a shortest cycle containing $e$ in $G$. 
Write $\ell_G(e)$ as $\ell(e)$ when $G$ is clear 
from the context.
Thus, the girth $g$ of $G$ is the minimum value of $\ell(e)$
among all edges $e$ in $G$.
Kaul and Mudrock~\cite{kaul1} showed that 
if $G$ has an even girth, then $\dpl{G}$.
We apply Theorem~\ref{th3-1} to 
generalize this result below.

\begin{theo}\label{th1-1}
For any graph $G$, 
if $\ell(e)$ is even for some edge $e$ in $G$,
then  $\dpl{G}$.
\end{theo} 
 
The converse statement of Theorem~\ref{th1-1} fails, and 
counterexamples will be given in Section~\ref{nsec4}.

\begin{theo}\label{th1-1-0}
There exist infinitely many graphs $G$  \label{cr9}
 such that   $\dpl{G}$ and 
$\ell(e)=3$ for each edge $e$ in $G$. 
\end{theo}

For a disconnected graph $G$,
if $\dpl{G_i}$ for some component $G_i$ of $G$,
then $\dpl{G}$ obviously holds.
This conclusion also holds for
connected graphs.

\begin{theo}\label{th1-1-00}
For a connected graph $G$, 
if $\dpl{G_i}$ for some block $G_i$ of $G$,
then $\dpl{G}$ holds.
\end{theo}

Some results on the study of graphs with the property $\dpg{G}$ 
have been obtained. 
Kaul and Mudrock \cite{kaul1} showed that 
$\dpg{G}$ holds for the graph $G$ 
obtained from any two odd cycle graphs 
$C_{2k+1}$ and $C_{2r+1}$ by 
identifying one edge in $C_{2k+1}$ 
with one edge in $C_{2r+1}$.
For two vertex-disjoint graphs $G$ and $G'$,
let $G\vee G'$ denote the join of $G$ and $G'$,
i.e., the graph obtained from $G$ and $G'$ 
by adding all edges in $\{uv: u\in V(G),v\in V(G')\}$.
Kaul and Mudrock \cite{kaul1} asked that 
for every graph $G$, does there exist $p\in \N$ 
such that $P_{DP}(K_p\vee G)\approx P(K_p\vee G)$,
where $K_p$ is the complete graph with $p$ vertices?
Recently, Mudrock and Thomason~\cite{mudrock1} 
showed that the problem has a positive answer for $p=1$.
Obviously, a graph is isomorphic to $K_1\vee G$ 
for some graph $G$ 
if and only if it has a dominating vertex
(i.e., a vertex which is adjacent to all other vertices in the graph).

For any graph $G$ and any integer $m>0$,
there is a special $m$-fold cover of $G$ 
which corresponds to proper $m$-colorings.
Let $\seth_0(G, m)$ denote the 
$m$-fold cover $(L,H)$ of $G$, where 
$L(u)=\{(u,i): i\in \brk{m}\}$ for each $u\in V(G)$ 
and $E_H(L(u),L(v))=\{(u,i)(v,i): i\in \brk{m}\}$
for each edge $uv$ in $G$.
The graph $H$ in $\seth_0(G, m)=(L,H)$
is denoted by $H_0(G,m)$ (or simply $H_0(m)$).
Obviously, $P_{DP}(G,\seth_0(G, m))=P(G,m)$ holds 
for each $m\in \N$.

Let $\DPg$ denote the set of graphs $G$ 
for which there exists $M\in \N$
such that for every $m$-fold cover $\seth=(L,H)$ of $G$, 
if $H\not\cong H_0(G,m)$, then 
$P_{DP}(G,\seth)>P(G,m)$ holds for all integers $m\ge M$.
By definition,  
$\dpg{G}$ holds for each graph $G$ in $\DPg$.
But it is unknown if the converse statement is also true.

\prom{prom01}
{Is it true that if $\dpg{G}$, then 
$G\in \DPg$?
}   

Our next result provides a sufficient condition 
for a graph $G$ to be in $\DPg$ 
and therefore $\dpg{G}$ holds.

\begin{theo}\label{th1-2}
If a graph $G$ contains 
a spanning tree $T$ such that 
for each edge $e$ in $E(G)\setminus E(T)$, 
$\ell(e)$ is odd and $e$ is 
contained in a cycle $C$ 
of length $\ell (e)$ with the property that 
$\ell(e')<\ell(e)$ holds for each 
$e'\in E(C)\setminus (E(T)\cup \{e\})$,
then $G\in \DPg$ and hence $\dpg{G}$.
\end{theo}

A vertex $u$ in a graph $G$ is called {\it simplicial}
if either $d_G(u)=0$ or $G[N(u)]$ is a complete graph.
A graph $G$ is called {\it chordal} if for each cycle $C$
in $G$, $G[V(C)]$ contains $3$-cycles.
Due to Dirac~\cite{dirac}, a graph $G$ 
	is chordal if and only if  
	there exists an ordering $v_1,v_2,\cdots,v_n$ of its vertices,
	called a {\it perfect elimination ordering},  
	such that each $v_i$ is simplicial in the subgraph of $G$ induced by $\{v_j: j\in \brk{i}\}$.
Due to Kaul and Mudrock \cite{kaul1},  \label{cr10}
for any chordal graph $G$,  
$P_{DP}(G,m)=P(G,m)$ holds for all $m\in \N$, 
and hence $\dpg{G}$ holds. 
We notice that this conclusion does not follow from
Theorem~\ref{th1-2}.
But the next result is its generalization. 

\begin{theo}\label{th1-2-0}
For any graph $G$ with a simplicial vertex $u$,
if $\dpg{G-u}$, then $\dpg{G}$; also,    \label{cr11}
if $G-u\in \DPg$, then $G\in \DPg$.
\end{theo}

Theorems~\ref{th1-1-0} and~\ref{th1-1-00}
are proved in Section~\ref{nsec4},
while 
Theorems~\ref{th1-1},~\ref{th1-2}
and~\ref{th1-2-0}
are proved in Sections~\ref{sec4}, 
~\ref{sec5} and~\ref{nsec5} respectively.

\section{Proof of Theorem~\ref{th1-1}
\label{sec4}}

The following result due to Kaul and Mudrock~\cite{kaul1}
will be applied to study graphs $G$ 
with the property  $\dpl{G}$.

\begin{theo}[\cite{kaul1}]\label{th3-1}
Let $G$ be a graph with an edge $e$. 
If $m\ge 2$ and 
$
P(G-e,m)<\frac{m}{m-1}P(G,m),
$
then $P_{DP}(G,m)<P(G,m)$.
\end{theo}

In this section, we shall apply two fundamental properties
of \label{cr12}
the chromatic polynomial $P(G,x)$ of $G$.
The variable $x$ in $P(G,x)$
can be considered  a real number.  
By the inclusion-exclusion principle, 
it can be proved that 
\equ{eq4-0}
{
P(G,x)=\sum_{A\subseteq E(G)}(-1)^{|A|}x^{c(A)},
}
where $c_G(A)$ (or simply $c(A)$) is the number 
of components in the spanning subgraph $\spnn{G}{A}$ of $G$ (see \cite{Whi1932}).
Note that (\ref{eq4-0}) holds even if $G$ 
has parallel edges or loops.

The deletion-contraction theorem of chromatic polynomials 
(see \cite{dong0,rea1,rea2}) states that for each 
edge $e$ in a graph $G$,
\equ{eq4-1}
{
P(G,x)=P(G-e,x)-P(G\slash e,x),
}
where $G\slash e$ is the graph obtained by contracting edge 
$e$ (i.e., the graph obtained from $G-e$ 
by identifying the two ends of $e$).
Clearly, $G\slash e$ may have parallel edges.
By (\ref{eq4-1}), for any $e\in E(G)$, 
when $x\ne 1$,  \label{oc4}
\begin{eqnarray}\label{eq4-2}
P(G-e,x)-\frac{x}{x-1}P(G,x) 
&=&P(G-e,x)-
\frac{x}{x-1}
\left (P(G-e,x)-P(G\slash e,x)
\right )
\nonumber \\
&=&\frac{1}{x-1}
\left (xP(G\slash e,x)-P(G-e,x)
\right ).
\end{eqnarray}

For any edge $e$ in $G$, 
let $\setc(e)$ denote the set of cycles in $G$ 
that \label{cr14}
contain $e$ and are of length $\ell(e)$.
Obviously,  $\setc(e)\ne \emptyset$ if $e$ is not a bridge of $G$.

\prop{pro4-1}
{Let $G$ be a simple graph and $e$ be an edge in $G$ 
with $\ell(e)<\infty$.
Then, the leading term in the polynomial 
$xP(G\slash e,x)-P(G-e,x)$ is  
$
(-1)^{\ell(e)-1}|\setc(e)|x^{n-\ell(e)+2}.
$
}

\proof 
Note that 
$G-e$ and $G\slash e$ have the same edge set,
i.e., $E(G)\setminus \{e\}$, and when $\ell(e)=3$, 
$G\slash e$ has parallel edges.
Applying (\ref{eq4-0}) to both 
$G-e$ and $G\slash e$, we have 
\equ{eq4-6}
{
P(G-e,x)=\sum_{A\subseteq E(G)
\setminus \{e\}}(-1)^{|A|}x^{c_G(A)}
}
and
\equ{eq4-7}
{
P(G\slash e,x)=\sum_{A\subseteq E(G)\setminus \{e\}}
(-1)^{|A|}x^{c_{G\slash e}(A)}.
}
Let $u,v$ be the two ends of $e$,  
and let $\cone$ be the set of subsets 
$A$ of $E(G)\setminus \{e\}$,
such that $u$ and $v$ are in the same 
component of the spanning subgraph $\spnn{G}{A}$ of $G$.
Let $\cone'$ be the set of 
subsets $A$ of $E(G)\setminus \{e\}$
with $A\notin \cone$.
If $A\in \cone$, then $c_G(A)=c_{G\slash e}(A)$;
and if $A\in \cone'$, then $c_G(A)=c_{G\slash e}(A)+1$.
Thus, (\ref{eq4-6}) and (\ref{eq4-7}) imply that 
\equ
{eq4-8}
{
xP(G\slash e,x)-P(G- e,x)
=
\sum_{A\in \cone}
(-1)^{|A|}x^{c_{G}(A)}(x-1).
}
For each $A\in \cone$, 
let $G_A$ denote the component of $\spnn{G}{A}$ 
that contains both vertices $u$ and $v$.
Then $G_A$ has a $(u,v)$-path $P$,
implying that $|V(G_A)|\ge |V(P)|\ge \ell(e)$.
If $|V(G_A)|=\ell(e)$, 
then $V(G_A)=V(P)$ and 
$P$ must be a path $C-e$ for 
some cycle $C\in \setc(e)$.
As each cycle in $G$ containing $e$ must be of length 
at least $\ell(e)$, 
$|V(G_A)|=\ell(e)$ implies that 
$G_A$ is a path $C-e$ for some cycle $C\in \setc(e)$.

Consequently, for each $A\in \cone$, 
$c_G(A)\le n-\ell(e)+1$ holds,
and $c_G(A)=n-\ell(e)+1$ 
if and only if $A\cup \{e\}$ is the edge set of 
some cycle $C$ in  $\setc(e)$.
Thus, 
$c_G(A)=n-\ell(e)+1$ holds for exactly
$|\setc(e)|$ subsets $A\in \cone$, 
and for each of them, $|A|=\ell(e)-1$.

By (\ref{eq4-8}) and the above conclusions, 
$xP(G\slash e,x)-P(G- e,x)$ 
is a polynomial of degree $n-\ell(e)+2$ 
and the coefficient of its leading term
is $(-1)^{\ell(e)-1}|\setc(e)|$.

Hence the result holds.
\proofend

We are now going to prove Theorem~\ref{th1-1}.

\noindent {\it Proof of Theorem~\ref{th1-1}}:
Let $e$ be an edge in $G$ such that $\ell(e)$ is even.
By the equality of (\ref{eq4-2})
and Proposition~\ref{pro4-1}, 
there exists $M\in \N$ such that
$P(G-e,m)<\frac{m}{m-1}P(G,m)$ for all integer $m\ge M$.
The result then follows 
from Theorem~\ref{th3-1}.
\proofend

\section{Proof of Theorems~\ref{th1-1-0} and~\ref{th1-1-00}
\label{nsec4}}

Let $\omega(G)$ denote the {\it clique number} of a graph $G$.
For any vertex-disjoint graphs $G_1$ and $G_2$
and $k\in \N$ with $k\le \min\{\omega(G_i):i=1,2\}$,
let ${\mathscr G}(G_1\cup_k G_2)$ denote the set of 
graphs obtained from $G_1$ and $G_2$ 
by identifying a $k$-clique in $G_1$ with 
a $k$-clique in $G_2$.
Due to Zykov~\cite{Zykov},  the following identity 
on chromatic polynomials holds for 
any $G\in {\mathscr G}(G_1\cup_k G_2)$ 
and all $m\ge k$:
\equ{zyk1}
{
P(G,m)=\frac{P(G_1,m)P(G_2,m)}{m(m-1)\cdots (m-k+1)}.
}
If $u$ is a simplicial vertex of a graph $G$, 
the following identity on chromatic polynomials
follows from (\ref{zyk1})
(also see \cite{dong0,rea2}): 
\equ{eq2-1}
{
P(G,m)=(m-d_G(u))P(G-u,m), \qquad \forall m\in \N.
}
It is natural to ask if (\ref{eq2-1}) 
holds for the DP color function. \label{cr15}

\prom{prob1}
{If $u$ is a simplicial vertex of $G$,
is it true that for 
all integers $m\ge d(u)$, \label{oc5}
\equ{eq2-2}
{
P_{DP}(G,m)=(m-d(u))P_{DP}(G-u,m)? 
}
}

As $P_{DP}(G,m)\ge (m-d(u))P_{DP}(G-u,m)$ by definition, 
to prove the equality of (\ref{eq2-2}), 
it suffices to show that 
$P_{DP}(G,m)\le (m-d(u))P_{DP}(G-u,m)$
for all integers $m\ge d(u)$. 
It is trivial that 
Problem~\ref{prob1} has a positive answer when $d(u)=0$,
and due to Theorem~\ref{beck-th},
it also has a positive answer when $d(u)=1$.
In this section, we show that it 
has a positive answer for $d(u)=2$.
Applying this conclusion, 
we are able to prove that 
the converse statement of Theorem~\ref{th1-1} fails.

\prop{pro2-2}
{
If $u$ is a simplicial vertex of $G$ with $d(u)=2$, then 
for each integer $m\ge 2$,
\equ{eq2-6}
{
P_{DP}(G,m)=(m-2)P_{DP}(G-u,m).
}
}

\proof 
Let $m\ge 2$.
If $m<\chi_{DP}(G-u)$, then 
$m<\chi_{DP}(G-u)\le \chi_{DP}(G)$,
implying that $P_{DP}(G-u,m)=P_{DP}(G,m)=0$.
It follows that (\ref{eq2-6}) holds in this case.

As $u$ is a simplicial vertex of $G$ with degree $2$, 
$\chi_{DP}(G)\ge \chi(G)\ge 3$,
implying that $P_{DP}(G,2)=0$.
Thus (\ref{eq2-6}) also holds when $m=2$.

Now let $m\ge \max\{3, \chi_{DP}(G-u)\}$ and 
let $\seth'=(L',H')$ be an $m$-fold cover of $G-u$ 
such that $P_{DP}(G-u,\seth')=P_{DP}(G-u,m)$ and 
$|E(H')|$ has the maximum value.  
It is clear that
$E_{H'}(L'(v_1),L'(v_2))$ is a  
matching in $H'$ of size $m$ 
for each pair of adjacent vertices $v_1$ and $v_2$ in $G-u$.

Let $N_G(u)=\{u_1,u_2\}$. 
Assume that $(u_1,j)$ and $(u_2,\pi(j))$ are
adjacent in $H'$ for each $j\in \brk{m}$,
where $\pi$ is a bijection from 
$\brk{m}$ to $\brk{m}$.

Let $H$ be the graph obtained from $H'$ 
and a complete graph 
with vertex set $\{(u,j):j\in \brk{m}\}$
by adding edges joining $(u,j)$ to 
both $(u_1,j)$ and $(u_2,\pi(j))$ 
for each $j\in \brk{m}$.
Let $\seth=(L,H)$ be the $m$-fold cover of $G$,
where $L(u)=\{(u,j):j\in \brk{m}\}$ 
and $L(v)=L'(v)$ for all $v\in V(G)-\{u\}$.

Let $I$ be any member in $\seti(H')$.
Assume that $(u_1,j_1)\in I\cap L(u_1)$
and $(u_2,\pi(j_2))\in I\cap L(u_2)$.
As  $(u_1,j_1)$ and $(u_2,\pi(j_1))$ are adjacent in $H$,
$j_1\ne j_2$.
Then, $I$ can be extended to exactly $(m-2)$ 
independent sets of $H$ of the form $I\cup \{(u,j)\}$,
where $j\in \brk{m}\setminus \{j_1,j_2\}$.
Thus,
\equ{eq2-7}
{
P_{DP}(G,\seth)=(m-2)P_{DP}(G-u,\seth')
=(m-2)P_{DP}(G-u,m),
}
by which $P_{DP}(G,m)\le (m-2)P_{DP}(G-u,m)$.
On the other hand, it is obvious that 
$P_{DP}(G,m)\ge (m-2)P_{DP}(G-u,m)$.
Thus, the result follows.
\proofend

\begin{figure}[h]
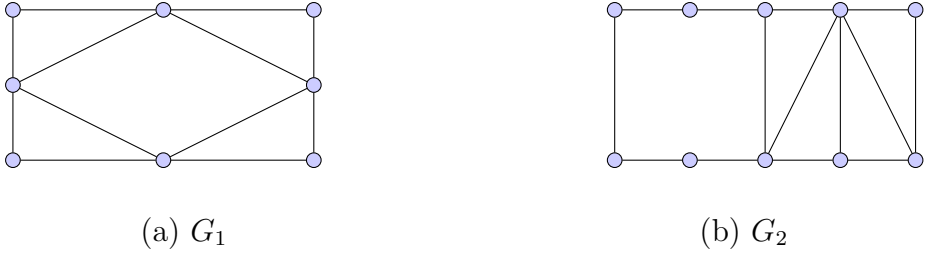


\tikzp{1}
{
\foreach \place/\z in 
{{(-6,0)/1}, {(-4,0)/2},{(-2,0)/3},
{(-6,1)/4},{(-2,1)/5},
{(-6,2)/6}, {(-4,2)/7},{(-2,2)/8}} 
\node[cblue] (a\z) at \place {};

\draw[black] (a1) -- (a2) -- (a3) -- (a5) -- (a8) -- (a7) 
-- (a6) -- (a4) -- (a1);

\draw[black] (a2) -- (a5) -- (a7) -- (a4) -- (a2);

\foreach \place/\z in 
{{(2,0)/1},{(3,0)/2},{(4,0)/3}, {(5,0)/4}, {(6,0)/5}, {(2,2)/6}, {(3,2)/7}, {(4,2)/8}, {(5,2)/9}, {(6,2)/10}}  
\node[cblue] (b\z) at \place {};

\draw[black] (b1) -- (b2) -- (b3) -- (b4) -- (b5) -- (b10) -- (b9) -- (b8) -- (b7) -- (b6) -- (b1); 
\path[black] (b9) edge (b3) edge (b4) edge (b5);

\draw[black] (b3) -- (b8);
}

\centerline{(a) $G_1$ \hspace{6 cm} 
(b) $G_2$}

\caption{$G_1\in \Phi(C_4)$ and $G_2\in \Phi(C_6)$}
\label{f4}
\end{figure}

For any graph $Q$ with at least one edge, 
let $\Phi(Q)$ be the family of graphs defined below:
\begin{enumerate}
\item   $Q\in \Phi(Q)$;
and 

\item  if $Q'\in \Phi(Q)$, 
then ${\mathscr G}(Q'\cup_2 K_3)\subseteq 
\Phi(Q)$. 
\end{enumerate}
For example, $G_1\in \Phi(C_4)$ and $G_2\in \Phi(C_6)$,
where $C_k$ is the cycle graph
of length $k$, and $G_1$ and $G_2$ are graphs in 
Figure~\ref{f4}.

By (\ref{eq2-1}) and Proposition~\ref{pro2-2}, 
for any graph $G\in \Phi(Q)$ and any integer $m\ge 2$, 
\equ{eq4-20}
{
P(G,m)=(m-2)^{|V(G)|-|V(Q)|}P(Q,m),
\quad P_{DP}(G,m)=(m-2)^{|V(G)|-|V(Q)|}P_{DP}(Q,m).
}
\label{cr18}

By (\ref{eq4-20}), we have the following observation.

\prop{pro2-3}
{
For any graph $Q$ with at least one edge and any $G\in \Phi(Q)$, 
if $\dpg{Q}$, then $\dpg{G}$; also,
if $\dpl{Q}$, then $\dpl{G}$.
}

We can now easily prove Theorem~\ref{th1-1-0}.

\noindent {\it Proof of Theorem~\ref{th1-1-0}}:
Let $Q$ be any graph with $\dpl{Q}$. 
Clearly, $Q$ contains edges. 
By Proposition~\ref{pro2-3},
$\dpl{G}$ holds for every $G\in \Phi(Q)$.
By the definition of $\Phi(Q)$, 
there are infinitely \label{cr19}
many graphs $G\in \Phi(Q)$ 
such that $\ell_G(e)=3$ for each edge $e$ in $G$.
For example, for graph $G_1$ in Figure~\ref{f4} (a),
if $G$ is a graph in $\Phi(G_1)$, then   
$\dpl{G}$ and 
$\ell_G(e)=3$ holds for each edge $e$ in $G$.
 
Theorem~\ref{th1-1-0} 
holds.
\proofend

Theorem~\ref{th1-1-00} will be proved directly 
by applying the following result due to 
Becker, Hewitt, Kaul, Maxfield, Mudrock,
Spivey, Thomason and Wagstrom~\cite{Beck}.

\begin{theo}[\cite{Beck}]\label{beck-th}
For any connected graph $G$ with blocks 
$G_1,G_2, \cdots,G_r$, where $r\ge 2$, 
\equ{pro1-e1}
{
P_{DP}(G,m)\le \frac 1{m^{r-1}}
\prod_{i=1}^r P_{DP}(G_i,m).
}
\end{theo}
\label{oc6}

\noindent {\it Proof of Theorem~\ref{th1-1-00}}:
Let $G_1,G_2,\cdots, G_r$ be the blocks of $G$.
By identity (\ref{zyk1}) and 
Theorem~\ref{beck-th}, we have 
\equ{ns3-eq10}
{
P_{DP}(G,m)\le \frac{1}{m^{r-1}}
\prod_{i=1}^r P_{DP}(G_i,m)
\le \frac{1}{m^{r-1}}
\prod_{i=1}^r P(G_i,m)
=P(G,m).
}
By (\ref{ns3-eq10}), if $\dpl{G_i}$ for some $i$,
then $\dpl{G}$ holds.
\proofend

It is natural to ask the following problem.

\prom{prom10}
{
For a connected graph $G$, 
if $\dpg{G_i}$ holds for each block $G_i$ of $G$,
is it true that $\dpg{G}$?
}

\section{Proof of Theorem~\ref{th1-2}
\label{sec5}}

\subsection{A set of ordered pairs $(G,T)$, where $T$ is a spanning tree
of $G$
\label{sec5-1}}

Let $\gtp$ be the set of ordered pairs $(G,T)$, where 
$G$ is a connected graph and $T$ is a spanning tree of $G$ 
such that for each edge $e$ in $E(G)\setminus E(T)$,
$\ell_G(e)$ is odd and $e$ is contained in a cycle $C$ 
of length $\ell_G(e)$ with the property 
that $\ell_G(e')<\ell_G(e)$ 
holds for each $e'\in E(C)\setminus (E(T)\cup \{e\})$.

Note that $\gtp$ contains a subfamily 
$\gtp_0$ of ordered pairs $(G,T)$, 
where $T$ is a spanning tree of $G$ such that 
for each $e\in E(G)\setminus E(T)$,
$\ell(e)$ is odd and the fundamental cycle $C_T(e)$ of $e$
with respect to $T$  
is of length $\ell(e)$.

Let $\gt$ (resp. $\gt_0$) be 
the set of graphs $G$ such that 
$(G,T)\in \gtp$ (resp. $(G,T)\in \gtp_0$)
for some spanning tree $T$ of $G$.
For example, for $i=1,2$, $G_i\in \gt_0$,
where $G_1$ and $G_2$ (i.e., the Petersen graph) are 
the graphs\label{cr20} in Figure~\ref{f2-1} (a) and (b)
respectively.
It is also obvious that $\gt_0$ 
contains every graph that\label{cr21} 
has a dominating vertex.
But, it can be verified that 
$G_3$ in Figure~\ref{f2-1} (c) 
belongs to $\gt\setminus \gt_0$.

\begin{figure}[h]
\tikzp{1}
{
\foreach \place/\z in {{(-6.4,0)/1}, {(-5,0)/2},{(-3.6,0)/3},
{(-5,1)/4},{(-3.6,1)/5},
{(-6.4,2)/6}, {(-5,2)/7},{(-3.6,2)/8}}  
\node[cblue] (w\z) at \place {};

\draw[ultra thick, black] (w1) -- (w2) -- (w4) -- (w7) -- (w6);
\path[ultra thick, black] (w4) edge (w3) edge (w5) edge (w8); 

\draw[thin, black] (w1) -- (w6);

\draw[thin, black] (w2) -- (w3) -- (w5) -- (w8) -- (w7);

\filldraw[black] (-5,-0.25) circle (0pt) node[anchor=north] {(a) $G_1$};

\foreach \place/\x in {{(-1,0)/1}, {(1,0)/2},
 {(-0.7,0.7)/8},{(0.7,0.7)/10},
{(1.7,1)/3}, {(-1.7,1)/6}, {(0,1)/9},
{(0,1.75)/7},
{(1,2)/4}, {(-1,2)/5}}   
\node[cblue] (v\x) at \place {};

\path[ultra thick, black] (v9) edge (v7) edge (v8) edge (v10);   

\draw[ultra thick, black] (v5) -- (v8);   
\draw[ultra thick, black] (v4) -- (v10); 
\draw[ultra thick, black] (v2) -- (v3) --(v7)--(v6)--(v1); 

\draw[thin, black] (v10) -- (v1) -- (v2) -- (v8); 
\draw[thin, black] (v3) -- (v4) -- (v5) -- (v6); 

\filldraw[black] (0,-0.25) circle (0pt) node[anchor=north] {(b) $G_2$};


\foreach \place/\y in {{(4,0)/1}, {(5.2,0)/2},  {(6.4,0)/3}, 
{(3.4,0.7)/4},{(4.6,0.7)/5}, {(5.8,0.7)/6}, {(7,0.7)/7}, 
{(4,1.4)/8},{(5.2,1.4)/9},{(6.4,1.4)/10}, 
{(4.6,2.1)/11}, {(5.8,2.1)/12}}  
\node[cblue] (u\y) at \place {};

\draw[ultra thick, black] (u1) -- (u2) --(u3)--(u6)--(u9)--(u5); 
\draw[ultra thick, black] (u4) -- (u8) --(u11)--(u9)--(u12)--(u10)--(u7); 

\draw[thin, black] (u4) -- (u5) -- (u6) -- (u7); 
\draw[thin, black] (u1) -- (u5);

\filldraw[black] (5.2,-0.25) circle (0pt) node[anchor=north] {(c) $G_3$};
}
\caption{$G_i\in \gt_0\subset \gt$ for $i=1,2$ and $G_3\in \gt\setminus \gt_0$}
\label{f2-1}
\end{figure}

Let $(G,T)\in \gtp$.
By definition, $\ell_G(e)$ is odd for each $e\in E(G)\setminus E(T)$.
But, it does not guarantee directly 
that $\ell_G(e)$ is not even for any $e\in E(T)$. 
By Theorem~\ref{th1-1}, 
 if $\ell_G(e)$ is even for some $e\in E(T)$, then Theorem~\ref{th1-2} fails. 
Thus,  before 
proving Theorem~\ref{th1-2}, 
it is necessary to show that 
$\ell_G(e)$ is not even for every $e\in E(T)$.

For $(G,T)\in \gtp$, if $E(G)=E(T)$, 
let $\ell(G,T)=\infty$;
otherwise, let 
$\ell(G,T)=\max\limits_{e\in E(G)\setminus E(T)}\ell_G(e)$.

\begin{pro}\label{pro5-1}
Let $(G,T)\in \gtp$.
For each edge $e\in E(T)$, 
if $e$ is not a bridge of $G$, then  
$\ell_G(e)$ is odd 
and $\ell_G(e)\le \ell(G,T)$.
\end{pro} 

\proof We prove the result by induction on $|E(G)|$.
Note that $|E(G)|\ge |E(T)|$. 
The result is obvious when $|E(G)|\le |E(T)|+1$.
Now assume that $|E(G)|\ge |E(T)|+2$
and the result holds for every ordered pair 
$(G',T')\in \gtp$ with $|E(G')|\le |E(G)|-1$.

Choose an edge $e_1$ in $E(G)\setminus E(T)$
such that $\ell_G(e_1)=\ell(G,T)<\infty$.
Clearly, $T$ is a spanning tree of $G-e_1$.
We first show that $(G-e_1,T)\in \gtp$.

Let $G'$ denote $G-e_1$ and 
let $e$ be any edge in $E(G')\setminus E(T)$.
As $e\in E(G)\setminus E(T)$, 
by definition, $\ell_G(e)$ is odd and 
$e$ is contained in a cycle $C$ in $G$ of length $\ell_G(e)$ 
such that $\ell_G(e')<\ell_G(e)$ 
for each $e'\in E(C)\setminus (E(T)\cup \{e\})$. 
By the choice of $e_1$, $\ell_G(e_1)\ge \ell_G(e)$, 
implying that $e_1\notin E(C)$.
Thus, $C$ is in $G'$ and $\ell_G(e)=\ell_{G'}(e)$.

Hence, by definition, $(G',T)\in \gtp$
and $\ell_{G'}(e)=\ell_G(e)$ for each 
$e\in E(G')\setminus E(T)$,
implying that $\ell(G',T)\le \ell(G,T)$.

By inductive assumption, the conclusion holds 
for $(G',T)\in \gtp$.
Now suppose $e_0\in E(T)$ and \label{cr22} 
$e_0$ is not a bridge in $G$.
Then, either $e_0$ is a bridge of $G'$ 
or $\ell_{G'}(e_0)$ is odd.
Furthermore, if $e_0$ is not a bridge of $G'$,
then $\ell_{G'}(e_0)\le \ell(G',T)
\le \ell(G,T)=\ell_G(e_1)$.
We shall show that $\ell_G(e_0)$ is odd 
and $\ell_G(e_0)\le \ell(G,T)$.

\noindent {\bf Case 1}: 
$e_0$ is a bridge of $G'$ (i.e., $G-e_1$).

In this case, for each cycle $C$ in $G$,
either $E(C)\cap \{e_0,e_1\}=\emptyset$ 
or $\{e_0,e_1\}\subseteq E(C)$,
implying that $\ell_G(e_0)=\ell_G(e_1)=\ell(G,T)$ is odd.

\noindent {\bf Case 2}: $\ell_{G'}(e_0)$ is odd.

In this case, $\ell_{G'}(e_0)\le \ell(G',T)
\le \ell(G,T)=\ell_G(e_1)$.
If $\ell_G(e_0)<\ell_{G'}(e_0)$, then 
$e_0$ is contained in a cycle $C$ in $G$ with $|E(C)|=\ell_{G}(e_0)$.
Since $|E(C)|=\ell_G(e_0)<\ell_{G'}(e_0)$,
$C$ is not in $G'$ and thus $e_1\in E(C)$,
implying that $\ell_G(e_1)\le |E(C)|$.
Hence 
$$
\ell_G(e_1)\le |E(C)|=\ell_G(e_0)<\ell_{G'}(e_0)
\le \ell(G',T)
\le \ell_G(e_1),
$$
a contradiction.
Hence $\ell_G(e_0)=\ell_{G'}(e_0)$ is odd.
Obviously, $\ell_G(e_0)=
\ell_{G'}(e_0)\le \ell(G',T)\le \ell(G,T)$.

Hence the result holds.
\proofend

\noindent {\bf Remark}: 
From the proof of Proposition~\ref{pro5-1}, 
for any $(G,T)\in \gtp$, 
$G$ can be obtained from $T$ by adding a sequence of edges.
Actually, $G$ is the last graph $G_k$ in a sequence of graphs $G_0,G_1,G_2,\cdots,G_k$,
where $k=|E(G)|-|V(G)|+1$, $G_0=T$
and each graph $G_{i+1}$, where $0\le i\le k-1$, 
can be obtained from $G_i$ 
by adding a new edge joining 
two nonadjacent vertices $u$ and $v$ in $G_i$
in which there is a shortest 
$(u, v)$-path $P$ such that
$|E(P)| \ge \ell_{G_i} (e) - 1$ for each $e \in E(G_i) \setminus E(T)$ and 
$|E(P)| >\ell_{G_i} (e) - 1$ for each
$e \in E(P) \setminus E(T)$.

\subsection{Proof of Theorem~\ref{th1-2}
\label{sec5-2}}

We are now going to prove Theorem~\ref{th1-2}.

\setcounter{countclaim}{0}

\noindent {\it Proof of Theorem~\ref{th1-2}}: 
Let $G\in \gt$ and $n=|V(G)|$. 
The result is trivial for $n=1$.
Now assume that $n\ge 2$.
By definition,
$(G,T)\in \gtp$ for some spanning tree $T$ of $G$.
Thus, for each $e\in E(G)\setminus E(T)$,
$\ell_G(e)$ is odd and $e$ is contained in a cycle $C$
of length $\ell_G(e)$ 
with the property that 
$\ell_G(e')<\ell_G(e)$ holds for each 
$e'\in E(C)\setminus (E(T)\cup \{e\})$.

Let $\seth=(L,H)$ be any $m$-fold cover of $G$
such that $H\not\cong H_0(G,m)$.
As $T$ is a spanning tree of $G$, 
by Proposition 21 in \cite{kaul1},
we may assume that $L(v)=\{(v,j): j\in \brk{m}\}$ 
for each $v\in V(G)$
and $E_H(L(u), L(v))\subseteq \{(u,j)(v,j):j\in \brk{m}\}$ 
for each $uv\in E(T)$.
Note that relabeling vertices in $L(u)$ for any $u\in V(G)$
does not affect the condition that 
$H\not\cong H_0(G,m)$.

If $E_H(L(u), L(v))\subseteq \{(u,j)(v,j):j\in \brk{m}\}$ 
holds for each $uv\in E(G)\setminus E(T)$, 
then $H\not\cong H_0(G,m)$ implies that 
$H$ is a proper spanning subgraph of $H_0(G,m)$.
Without loss of generality, assume that 
$(u,1)(v,1)\notin E_H(L(u), L(v))$ for some edge 
$uv\in E(G)$.
Then, for $m\ge n-2$,
$$
P_{DP}(G, \seth)-P(G,m)
\ge P_{DP}(G-\{u,v\},\seth')>0, 
$$
where $\seth'=(L',H')$ is the $(m-1)$-fold cover of $G-\{u,v\}$, 
$L'(w)=L(w)\setminus \{(w,1)\}$ 
for all $w\in V(G)\setminus \{u,v\}$
and $H'=H[\cup_{w\in V(G)\setminus \{u,v\}}L'(w)]$.
Thus, the result holds in this case.
 
Now assume that 
$E_H(L(u), L(v))\not \subseteq \{(u,j)(v,j):j\in \brk{m}\}$ 
for some $uv\in E(G)\setminus E(T)$.
By definition, adding any possible edge to $H$ does not 
increase the value of $P_{DP}(G,\seth)$.
Thus, we can assume that $|E_H(L(u),L(v))|=m$
for each edge $uv\in E(G)$,
and in particular, 
$E_H(L(u), L(v))=\{(u,j)(v,j):j\in \brk{m}\}$ 
for each $uv\in E(T)$.

For each $e=uv\in E(G)$, 
let $X_e=E_H(L(u),L(v))\setminus 
\{(u,j)(v,j):j\in \brk{m}\}$. 
As $|E_H(L(u),L(v))|=m$,  
$X_e=\emptyset$ if and only if
$E_H(L(u),L(v))= \{(u,j)(v,j):j\in \brk{m}\}$.
By the assumption above, 
$X_e=\emptyset$ for each $e\in E(T)$,
but $X_e\ne \emptyset$
for some edge $e\in E(G)\setminus E(T)$.
For  $s\ge 3$, let 
\equ{eq5-3}
{
\scrx_s=\bigcup_{e\in E(G)\setminus E(T)\atop \ell_G(e)=s} X_e.
}
By the given condition, $\ell_G(e)\ge 3$ 
is odd for each $e\in E(G)\setminus E(T)$,
implying that $\scrx_s=\emptyset$ 
for each even $s\ge 4$.
Now assume that $r$ is the minimum integer such that 
$\scrx_r\ne \emptyset$. So $r\ge 3$ and $r$ is odd.
We will prove Theorem~\ref{th1-2} \label{cr23}
by an approach similar to the proof
of Theorem 7 in \cite{mudrock1}. 

We first find an expression for $P_{DP}(G,\seth)$ 
which is similar to (\ref{eq4-0})
for $P(G,m)$.
Let $\scrs$ be the set of subsets $S$ of $V(H)$
with $|S\cap L(v)|=1$ for each $v\in V(G)$.
For each edge $e=uv\in E(G)$, 
let $\scrs_e$ be the set of $S\in \scrs$
such that the two vertices in $S\cap (L(u)\cup L(v))$ are
adjacent in $H$.
For each $A\subseteq E(G)$, 
let 
\equ{eq5-4}
{
\scrs_A=\bigcap_{e\in A} \scrs_e.
} 
As $P_{DP}(G,\seth)=|\scrs|-|\cup_{e\in E(G)}\scrs_e|$, 
by the inclusion-exclusion principle, we have
\equ{eq5-5}
{
P_{DP}(G,\seth)=\sum_{A\subseteq E(G)}(-1)^{|A|} |\scrs_A|.
}

For each $U\subseteq V(G)$, let $\scrs|_U$ be the set of 
subsets $S$ of $V(H)$ such that 
$|S\cap L(v)|=1$ for each $v\in U$ 
and $S\cap L(v)=\emptyset$ for each $v\in V(G)\setminus U$.
For any subgraph $G_0$ of $G$ and $S\in \scrs|_{V(G_0)}$,
let $\spng{S}{G_0}$  denote 
the spanning subgraph of $H[S]$ 
with edge set $\{(u,j_1)(v,j_2)\in E(H):
uv\in E(G_0),u,v\in V(G_0),(u,j_1),(v,j_2)\in S\}$.
Equivalently, $\spng{S}{G_0}$ can be obtained from $H[S]$
by deleting all those edges $(u,j_1)(v,j_2)$ in $H[S]$
with $uv\notin E(G_0)$.\label{cr24}
Clearly, $\spng{S}{G_0}$ is $H[S]$
when $G_0$ is a subgraph of $G$ induced by $V(G_0)$.
For any $S\in \scrs|_{V(G_0)}$,
$|E(\spng{S}{G_0})|\le |E(G_0)|$ holds,
and the following statements are equivalent:

(a) $\spng{S}{G_0}\cong G_0$;

(b) $|E(\spng{S}{G_0})|= |E(G_0)|$; and 

(c) for each $uv\in E(G_0)$, 
the two vertices in $S\cap (L(u)\cup L(v))$ are
adjacent in $H$.

Let $\scrh(G_0)$ be the set of 
subgraphs $\spng{S}{G_0}$ of $H$, 
where $S\in \scrs|_{V(G_0)}$,
such that $\spng{S}{G_0}\cong G_0$.

Recall that for $A\subseteq E(G)$,
$\spnn{G}{A}$ is the spanning 
subgraph of $G$ with edge set $A$,
and $c(A)$ is the number of 
components of $\spnn{G}{A}$.
By the definition of $\scrs_A$,
the following claim holds.

\begin{clm}\label{clone}
For any $A\subseteq E(G)$, 
if $G_1,G_2,\cdots,G_{c(A)}$ are the components of 
$\spnn{G}{A}$,
then 
$$
|\scrs_A|=\prod_{i=1}^{c(A)} |\scrh(G_i)|.
$$
\end{clm} 

\begin{clm}\label{cltwo}
Let $G_0$ be a connected subgraph of $G$.
If $\spng{S_1}{G_0}, \spng{S_2}{G_0}
\in {\mathscr H}(G_0)$, 
where $S_1,S_2\in \scrs|_{V(G_0)}$,
then either $S_1=S_2$ or $S_1\cap S_2=\emptyset$.
Hence $|{\mathscr H}(G_0)|\le m$,
where the equality holds 
if $X_e=\emptyset$ holds for each edge $e\in E(G_0)$.
\end{clm}

\proofn
Suppose that $\spng{S_1}{G_0}, \spng{S_2}{G_0}
\in {\mathscr H}(G_0)$.
Then, $\spng{S_1}{G_0}\cong \spng{S_2}{G_0}\cong G_0$,
implying that whenever $uv\in E(G_0)$,
the two vertices in $S_i\cap (L(u)\cup L(v))$ are
adjacent in $H$ for $i=1,2$. 
Let $uv$ be any edge in $G_0$. 
As $E_H(L(u),L(v))$ is a matching of $H$ of size $m$,
each vertex in $L(u)$ is only adjacent to 
one vertex in $L(v)$.
If $\spng{S_1}{G_0}$ and $\spng{S_2}{G_0}$ have 
a common vertex in $L(u)$, 
then $\spng{S_1}{G_0}$ and $\spng{S_2}{G_0}$ must have 
a common vertex in $L(v)$.
As $G_0$ is connected, we conclude that 
either $S_1\cap S_2=\emptyset$ or $S_1=S_2$.
Thus, $|{\mathscr H}(G_0)|\le m$ holds.

If $X_e=\emptyset$ holds for each edge $e\in E(G_0)$,
then $\spng{S_j}{G_0}\in {\mathscr H}(G_0)$ 
for each $j\in \brk{m}$, where 
$S_j=\{(u,j): u\in V(G_0)\}$.
Thus, $|{\mathscr H}(G_0)|=m$ and  
Claim \ref{cltwo} holds.
\claimend

\begin{clm}
\label{clthree} 
Let $G_0$ be a connected subgraph of $G$.
If $X_e=\emptyset$ holds for 
each $e\in E(G_0)$ that \label{cr26}
is not a bridge of $G_0$, then $|\scrh(G_0)|= m$.
\end{clm}
 
\proofn
Assume that $X_e=\emptyset$ holds for 
each $e\in E(G_0)$ that is not a bridge of $G_0$.
Let $B$ be any block of $G_0$.
If $B$ is trivial
(i.e., it consists of a bridge $e=uv$ of $G_0$ only),
then, it is clear that $\scrh(B)$ has exactly $m$ 
members which correspond to the $m$ edges in 
$E_H(L(u),L(v))$.
If $B$ is an non-trivial block of $G_0$,
we have $X_e=\emptyset$ for each $e\in E(B)$,
and 
$\scrh(B)$ has exactly $m$ members 
$\spng{S_j}{B}$ for $j\in \brk{m}$,
where $S_j=\{(v,j): v\in V(B)\}$.
Thus, $|\scrh(G_0)|=m$ if $G_0$ has only one block. 

Suppose that $G_0$ has at least two blocks
and $B_0$ is a block of $G_0$ 
which has only one vertex 
$u$ shared by other blocks of $G_0$.
Let $G'$ denote $G_0-(V(B_0)\setminus \{u\})$.
Assume that 
both $\scrh(G')$ and $\scrh(B_0)$ have exactly $m$ members.
Each member $\spng{S'}{G'}$ of $\scrh(G')$ 
can be extended to exactly one member of 
$\scrh(G_0)$ by combining $\spng{S'}{G'}$ 
with the member in $\scrh(B_0)$
which shares a vertex in $L(u)$ with $\spng{S'}{G'}$.
Hence $|\scrh(G_0)|=m$. 

Claim \ref{clthree} holds.
\claimend 
 
The next claim follows from 
Claims \ref{clone}, \ref{cltwo} and \ref{clthree} directly.

\begin{clm}\label{clfour}
For each $A\subseteq E(G)$, 
$|\scrs_A|\le m^{c(A)}$ holds.
If $X_e=\emptyset$ holds for 
each $e\in A$ that is not a bridge of $\spnn{G}{A}$, 
then $|\scrs_A|= m^{c(A)}$.
\end{clm}

By Claim~\ref{clfour}, the next claim follows.

\begin{clm}\label{clsix} 
For any $A\subseteq E(G)$, 
if $|A|$ is odd, 
$(-1)^{|A|} \left (|\scrs_A|-m^{c(A)}\right )
=m^{c(A)}-|\scrs_A|\ge 0$.
\end{clm}

Let $\scre$ be the set of subsets $A$ of $E(G)$ such that 
$X_e\ne \emptyset$ holds for 
some $e\in A$ that is not a bridge of $\spnn{G}{A}$.
Note that such an edge \label{cr27}
$e$ may be not unique. 
By (\ref{eq4-0}), (\ref{eq5-5}) and 
Claim~\ref{clfour}, we have  
\equ{eq5-6}
{
P_{DP}(G,\seth)-P(G,m)=
\sum_{A\in \scre}(-1)^{|A|} 
\left (|\scrs_A|-m^{c(A)}\right ).
}  

The following claim presents some properties 
of members in $\scre$.
\label{oc8}

\begin{clm}\label{clfive}
For each $A\in \scre$,
$\spnn{G}{A}$ has a component $G_1$
and an edge $e$ in some cycle of $G_1$ 
with $X_e\ne \emptyset$. 
Furthermore, $|V(G_1)|\ge r$ and $c(A)\le n-r+1$,
and $|A|=r$ whenever $c(A)=n-r+1$.
\end{clm}

\proofn
As $A\in \scre$, $\spnn{G}{A}$ has an edge $e$
that is not a bridge of $\spnn{G}{A}$ 
such that $X_e\ne \emptyset$.
Let $G_1$ be the component of $\spnn{G}{A}$ containing $e$.
As $X_e\ne \emptyset$, we have $\ell_G(e)\ge r$.
Thus, each cycle in $G_1$ containing $e$ has at least $r$ edges,
implying that $|V(G_1)|\ge r$, and hence $c(A)\le n-r+1$.

Assume that $c(A)=n-r+1$. 
Then $|V(G_1)|\ge r$ implies that $|V(G_1)|=r$
and all other components of $\spnn{G}{A}$
are isolated vertices.
As $e$ is in a cycle $C$ of length $r$ in $G_1$
and each cycle containing $e$ is of length at 
least $r$, $|V(G_1)|=r$ implies that 
$G_1\cong C$ and $|A|=r$.
 
Claim~\ref{clfive} holds.
\claimend 

Assume that $\{e_1,e_2,\cdots,e_s\}$ is the set of 
edges in $E(G)\setminus E(T)$ with $\ell_G(e_i)=r$.
By the given condition, 
for each $i\in \brk{s}$, 
$e_i$ is contained in a cycle, denoted by $C_i$,
such that $|V(C_i)|=r$ and $\ell_G(e')<r$ 
for each $e'\in E(C_i)\setminus (E(T)\cup \{e_i\})$.
Thus,  $E(C_i)\cap \{e_j: j\in \brk{s}\}=\{e_i\}$
for each $i\in \brk{s}$,
implying that $C_1,C_2,\cdots,C_s$ are pairwise distinct.

\begin{clm}
\label{clseven-1}
For each $i\in \brk{s}$, 
$\left |\scrh(C_i)\right |=m-|X_{e_i}|$.
\end{clm}

\proofn
Without loss of generality, 
let $V(C_i)=\{v_1,v_2,\cdots,v_r\}$
and let $v_1v_2\cdots v_r$  be 
the path $C_i-e_i$ in $G$.
Obviously, $e_i$ is the edge $v_1v_r$.
By the definition of $r$, \label{oc91}
$X_{e'}=\emptyset$ holds for 
each $e'\in E(G)\setminus E(T)$ 
with $\ell_G(e')<r$.
By the given condition on $C_i$, 
$X_{e'}=\emptyset$ holds for each edge $e'$ in the path $v_1v_2\cdots v_r$,
implying that the subgraph obtained from 
$H[S]$, where $S=\{(v_q,j): q\in [r],j\in \brk{m}\}$,
by removing all edges in $E_H(L(v_1),L(v_r))$,
consists of $m$ disjoint paths 
$(v_1,j)(v_2,j)\cdots (v_r,j)$
for $j=1,2,\cdots,m$.  
Assume that 
\equ{eq5-23}
{
E_H(L(v_1),L(v_r))\setminus X_{e_i}
=\{(v_1,j)(v_r,j): 1\le j\le m-|X_{e_i}|\}.
}
Then $(v_1,j)(v_r,j)\notin E(H)$ for each $j$ with  
$m-|X_{e_i}|<j\le m$.
Let $S_j=\{(v_q,j): q\in [r]\}$ for 
each $j\in \brk{m}$.
Clearly, $\spng{S_j}{C_i}\cong C_i$
if and only if $1\le j\le m-|X_{e_i}|$.
On the other hand, for any $S'\in \scrs|_{V(C_i)}$,
if $\spng{S'}{C_i}\cong C_i$,
then $\spng{S'}{C_i}$ must contain
a path $(v_1,j)(v_2,j)\cdots (v_r,j)$ for some 
$j\in \brk{m}$,
implying that $S'=S_j$ for some $j\in \brk{m}$.
Thus, $|\scrh(C_i)|=m-|X_{e_i}|$. 
\claimend 

Now we are going to apply 
Claims~\ref{clsix} and~\ref{clfive} to
prove the next claim.

\begin{clm}\label{clseven} 
The following result holds:
\equ{eq5-8}
{
\sum_{A\in \scre\atop c(A)=n-r+1}(-1)^{|A|} 
\left (|\scrs_A|-m^{c(A)}\right )
\ge |\scrx_r| m^{n-r}.
}
\end{clm}

\proofn
By Claim~\ref{clfive}, $|A|=r$ 
for each $A\in \scre$ with 
$c(A)=n-r+1$.
As $r$ is odd, by Claim~\ref{clsix}, 
for any  
$\scre_0\subseteq \{A\in \scre: c(A)=n-r+1\}$,  
we have 
\equ{eq5-11-0}
{
\sum_{A\in \scre\atop c(A)=n-r+1}
(-1)^{|A|} 
\left (|\scrs_A|-m^{c(A)}\right ) 
\ge \sum_{A\in \scre_0}
\left (m^{n-r+1}-|\scrs_A|\right ).
}

For each $i\in \brk{s}$,  
$\spnn{G}{E(C_i)}$ consists of 
exactly $n-r+1$ components,
i.e.,  $C_i$ and $n-r$  isolated vertices 
in $V(G)\setminus V(C_i)$.
By Claim~\ref{clseven-1},
$|\scrh(C_i)|=m-|X_{e_i}|$.
Thus, by Claim~\ref{clone},   
\equ{eq5-11}
{
|\scrs_{E(C_i)}|=|\scrh(C_i)| m^{n-r}
= (m-|X_{e_i}|)m^{n-r}=m^{n-r+1}-|X_{e_i}|m^{n-r}.
}
Let $\scre_0=\{E(C_i):i\in \brk{s}\}$.
By (\ref{eq5-11-0}) and (\ref{eq5-11}),
\equ{eq5-7}
{
\sum_{A\in \scre\atop c(A)=n-r+1}
(-1)^{|A|} 
\left (|\scrs_A|-m^{c(A)}\right ) 
\ge \sum_{i=1}^s
\left (m^{n-r+1}-|\scrs_{E(C_i)}|\right ) 
=  \sum_{i=1}^s
|X_{e_i}|m^{n-r}
=|\scrx_r| m^{n-r}.
}
Claim~\ref{clseven} holds.
\claimend

\begin{clm}\label{cleight-1}
For any subgraph $G_1$ of $G$,
if $\ell_G(e)\le r$ for each edge $e\in E(G_1)$,
then $|\scrh(G_1)|\ge m-2|\scrx_r|$.
\end{clm}

\proofn
For each $j\in \brk{m}$, 
let $S_j=\{(u,j): u\in V(G_1)\}$
and $Q_j=\spng{S_j}{G_1}$.
By the definition of $\spng{S_j}{G_1}$, 
$Q_j\in \scrh(G_1)$
if and only if 
$(u,j)(v,j)\in E(H)$  for each $uv\in E(G_1)$.

Let $S=\cup_{j\in \brk{m}}S_j$, and 
let $\psi: S
\rightarrow \{0,1\}$ 
be the mapping defined below:
\equ{eq5-21}
{
\psi((u,j))=
\left \{
\begin{array}{ll}
1, \qquad &\mbox{if }
(u,j)(v,j')\in E(H) \mbox{ for some }
v\in N_{G_1}(u) \mbox{ and } j'\ne j;\\
0, &\mbox{otherwise}.
\end{array}
\right.
}
If 
$\psi((u,j))=1$, by definition, $(u,j)$
is one end of some edge $(u,j)(v,j')$ of $X_e$, 
where $e=uv\in E(G_1)$. 
Thus, 
\begin{eqnarray}
\label{eq5-22}
\sum_{(u,j)\in S}\psi((u,j))
&\le & \sum_{e\in E(G_1)} 
\sum_{(u,j)(v,j')\in X_e}(\psi((u,j))+\psi((v,j')))
\nonumber \\
&=& 2\sum_{e\in E(G_1)}|X_e| \nonumber \\
&\le & 2|\scrx_r|,
\end{eqnarray}
where the last inequality follows from the 
facts that for each $e\in E(G_1)$, 
$\ell_G(e)\le r$ holds, 
and $\ell_G(e)<r$ implies that $X_e=\emptyset$.

By the definition of $\psi$, 
$Q_j\not \cong G_1$ if and only if 
$\psi((u,j))=1$ for some $u\in V(G_1)$.
Then, by (\ref{eq5-22}), there are at most $2|\scrx_r|$
numbers $j\in \brk{m}$ such that $Q_j\not \cong G_1$,
implying that 
\equ{eq5-12}
{
|\scrh(G_1)| 
\ge m-2|\scrx_r|.
}
Thus, Claim~\ref{cleight-1} holds.
\claimend

\begin{clm}\label{cleight}
For any $A\in \scre$ with $c(A)=n-r$,
we have  
$|\scrs_A|\ge (m-2|\scrx_r|)m^{m-r-1}$.
\end{clm}

\proofn
Let $A\in \scre$ with $c(A)=n-r$. 
By Claim~\ref{clfive}, $\spnn{G}{A}$ has a component $G_1$ 
with $|V(G_1)|\ge r$. 
Let $G_2,\cdots, G_{n-r}$ be the components of $\spnn{G}{A}$
different from $G_1$
with $|V(G_2)|\ge \cdots \ge |V(G_{n-r})|$. 
As $c(A)=n-r$, one of the two cases below happens:
\begin{enumerate}
\item $|V(G_1)|=r$, $|V(G_2)|=2$
and  $|V(G_i)|=1$ for all $3\le i\le n-r$,
or

\item $|V(G_1)|=r+1$ 
and $|V(G_i)|=1$ for all $2\le i\le n-r$.
\end{enumerate}

In both Cases (i) and (ii) above, by Claim~\ref{clthree}, 
$|\scrh(G_i)|=m$ holds for all $i=2,3,\cdots,n-r$.
By Claim~\ref{clone}, 
it remains to show that 
$|\scrh(G_1)|\ge m-2|\scrx_r|$.

\begin{figure}[h]
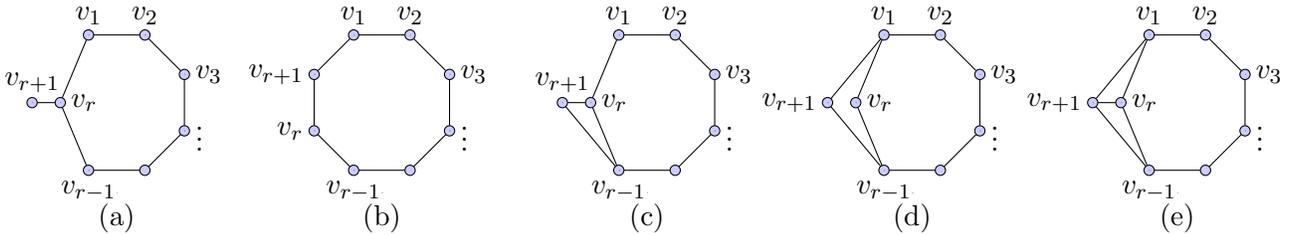


\tikzstyle{cblue}=[circle, draw, thin,fill=blue!20, scale=0.35]
\tikzp{0.75}
{
\foreach \place/\x in {{(6,0)/1}, {(7,0)/2},  
{(7.7,0.7)/3}, 
{(5,1.2)/4}, {(5.5,1.2)/5},
{(7.7,1.7)/6}, 
{(6,2.4)/7}, {(7,2.4)/8}}  
\node[cblue] (e\x) at \place {};

\draw[black] (e1) -- (e2) -- (e3) -- (e6) -- (e8) -- (e7) -- (e4) -- (e5) -- (e1);
\draw[black] (e1) -- (e4);
\draw[black] (e5) -- (e7);

\filldraw[black] (e1) circle (0pt) node[anchor=north] {\small $v_{r-1}$};
\filldraw[black] (e3) circle (0pt) node[anchor=west] {\small $\vdots$};
\filldraw[black] (e6) circle (0pt) node[anchor=west] {\small $v_3$};
\filldraw[black] (e7) circle (0pt) node[anchor=south] {\small $v_1$};
\filldraw[black] (e8) circle (0pt) node[anchor=south] {\small $v_2$};
\filldraw[black] (e5) circle (0pt) node[anchor=west] {\small $v_r$};
\filldraw[black] (e4) circle (0pt) node[anchor=east] {\small $v_{r+1}$};

\filldraw[black] (6.5,-0.4) circle (0pt) node[anchor=north] {\small (e)};

\foreach \place/\x in {{(1.3,0)/1}, {(2.3,0)/2},  
{(3,0.7)/3}, 
{(0.3,1.2)/4}, {(0.8,1.2)/5},
{(3,1.7)/6}, 
{(1.3,2.4)/7}, {(2.3,2.4)/8}}  
\node[cblue] (d\x) at \place {};

\draw[black] (d1) -- (d2) -- (d3) -- (d6) -- (d8) -- (d7) -- (d4);
\draw[black] (d5) -- (d1);
\draw[black] (d1) -- (d4);
\draw[black] (d5) -- (d7);

\filldraw[black] (d1) circle (0pt) node[anchor=north] {\small $v_{r-1}$};
\filldraw[black] (d3) circle (0pt) node[anchor=west] {\small $\vdots$};
\filldraw[black] (d6) circle (0pt) node[anchor=west] {\small $v_3$};
\filldraw[black] (d7) circle (0pt) node[anchor=south] {\small $v_1$};
\filldraw[black] (d8) circle (0pt) node[anchor=south] {\small $v_2$};
\filldraw[black] (d5) circle (0pt) node[anchor=west] {\small $v_r$};
\filldraw[black] (d4) circle (0pt) node[anchor=east] {\small $v_{r+1}$};

\filldraw[black] (1.8,-0.4) circle (0pt) node[anchor=north] {\small (d)};

\foreach \place/\x in {{(-3.4,0)/1}, {(-2.4,0)/2},  
{(-1.7,0.7)/3}, 
{(-4.4,1.2)/4}, {(-3.9,1.2)/5},
{(-1.7,1.7)/6}, 
{(-3.4,2.4)/7}, {(-2.4,2.4)/8}}  
\node[cblue] (c\x) at \place {};
 
\draw[black] (c1) -- (c2) -- (c3) -- (c6) -- (c8) -- (c7) --  (c5) -- (c1);
\draw[black] (c1) -- (c4) -- (c5);

\filldraw[black] (c1) circle (0pt) node[anchor=north] {\small $v_{r-1}$};
\filldraw[black] (c3) circle (0pt) node[anchor=west] {\small $\vdots$};
\filldraw[black] (c6) circle (0pt) node[anchor=west] {\small $v_3$};
\filldraw[black] (c7) circle (0pt) node[anchor=south] {\small $v_1$};
\filldraw[black] (c8) circle (0pt) node[anchor=south] {\small $v_2$};
\filldraw[black] (c5) circle (0pt) node[anchor=west] {\small $v_r$};
\filldraw[black] (c4) circle (0pt) node[anchor=south] {\small $v_{r+1}$};

\filldraw[black] (-2.9,-0.4) circle (0pt) node[anchor=north] {\small (c)};

\foreach \place/\x in {{(-8.1,0)/1}, {(-7.1,0)/2},  
{(-6.4,0.7)/3}, 
{(-8.8,1.7)/4}, {(-8.8,0.7)/5},
{(-6.4,1.7)/6}, 
{(-8.1,2.4)/7}, {(-7.1,2.4)/8}}  
\node[cblue] (b\x) at \place {};
 
\draw[black] (b1) -- (b2) -- (b3) -- (b6) -- (b8) -- (b7) -- (b4) -- (b5) -- (b1);

\filldraw[black] (b1) circle (0pt) node[anchor=north] {\small $v_{r-1}$};
\filldraw[black] (b3) circle (0pt) node[anchor=west] {\small $\vdots$};
\filldraw[black] (b6) circle (0pt) node[anchor=west] {\small $v_3$};
\filldraw[black] (b7) circle (0pt) node[anchor=south] {\small $v_1$};
\filldraw[black] (b8) circle (0pt) node[anchor=south] {\small $v_2$};
\filldraw[black] (b5) circle (0pt) node[anchor=east] {\small $v_r$};
\filldraw[black] (b4) circle (0pt) node[anchor=east] {\small $v_{r+1}$};

\filldraw[black] (-7.6,-0.4) circle (0pt) node[anchor=north] {\small (b)};

\foreach \place/\x in {{(-12.8,0)/1}, {(-11.8,0)/2},  
{(-11.1,0.7)/3}, 
{(-13.8,1.2)/4}, {(-13.3,1.2)/5},
{(-11.1,1.7)/6}, 
{(-12.8,2.4)/7}, {(-11.8,2.4)/8}}  
\node[cblue] (a\x) at \place {};

\draw[black] (a1) -- (a2) -- (a3) -- (a6) -- (a8) -- (a7) --  (a5) -- (a1);
\draw[black] (a4) -- (a5);

\filldraw[black] (a1) circle (0pt) node[anchor=north] {\small $v_{r-1}$};
\filldraw[black] (a3) circle (0pt) node[anchor=west] {\small $\vdots$};
\filldraw[black] (a6) circle (0pt) node[anchor=west] {\small $v_3$};
\filldraw[black] (a7) circle (0pt) node[anchor=south] {\small $v_1$};
\filldraw[black] (a8) circle (0pt) node[anchor=south] {\small $v_2$};
\filldraw[black] (a5) circle (0pt) node[anchor=west] {\small $v_r$};
\filldraw[black] (a4) circle (0pt) node[anchor=south] {\small $v_{r+1}$};

\filldraw[black] (-12.3,-0.4) circle (0pt) node[anchor=north] {\small (a)};

}
\caption{Possible structures of $G_1$ when $|V(G_1)|=r+1$}
\label{f3}
\end{figure}

In both cases above, by Claim~\ref{clfive},
there is an edge $e$ with $X_e\ne \emptyset$ 
which is in some cycle of $G_1$.
Such an edge may be not unique. 
As $X_e\ne \emptyset$, we have $\ell_G(e)\ge r$.
Thus, each cycle in $G_1$ containing $e$ 
must be of length at least $r$.
In Case (i),  
$G_1$ can only be a cycle of length $r$.
In Case (ii), 
it can be verified that 
$G_1$ is one of the graphs in Figure~\ref{f3}.

As each cycle in $G_1$ is of length at most $r+1$, 
for each edge $e'$ in cycles of $G_1$,
we have $\ell_G(e')\le r+1$.
As $\ell_G(e')$ is odd, we have 
$\ell_G(e')\le r$ for such edges $e'$.
Thus, if $G_1$ contains an edge 
$e'$ with $\ell_G(e')\ge r+2$,
then $e'$ must be a bridge of $G_1$.

If $G_1$ has no bridge, 
then $\ell_G(e)\le r$ for each $e\in E(G_1)$.
By Claim~\ref{cleight-1},
$\scrh(G_1)\ge m-2|\scrx_r|$.

If $G_1$ has bridges,
then $G_1$ is the graph in 
Figure~\ref{f3} (a), where $G_1-v_{r+1}$ is 
a cycle. 
By Claim~\ref{cleight-1} again,
$\scrh(G_1-v_{r+1})\ge m-2|\scrx_r|$.
Clearly, each member of $\scrh(G_1-v_{r+1})$
can be extended to a member of $\scrh(G_1)$,
even when $X_{v_rv_{r+1}}\ne \emptyset$. \label{cr28}
Thus, Claim~\ref{cleight} also holds in this case.

Claim~\ref{cleight} is proved. 
\claimend

For any $k\in \brk{n-r}$, 
let $\phi_k$ be the number of  
elements $A$ of $\scre$ such that $c(A)=k$ and $|A|$ is even.

\begin{clm}\label{clnine}
The following inequality holds: 
\equ{eq5-9}
{
\sum_{A\in \scre\atop c(A)=n-r}
(-1)^{|A|} 
\left (|\scrs_A|-m^{c(A)}\right )
\ge -2\phi_{n-r}  |\scrx_r|  m^{n-r-1}.
}
\end{clm}

\proofn
By Claim~\ref{clsix},
\equ{eq5-10}
{
\sum_{A\in \scre\atop c(A)=n-r}
(-1)^{|A|} \left (|\scrs_A|-m^{c(A)}\right )
\ge \sum_{A\in \scre, c(A)=n-r
\atop |A|\ \mbox{\small is even}} 
\left (|\scrs_A|-m^{c(A)}\right ).
}
For each $A\in \scre$ with $c(A)=n-r$, by Claim~\ref{cleight},
\equ{eq5-10-1}
{
|\scrs_A|-m^{c(A)}\ge 
(m-2|\scrx_r|)m^{n-r-1}-m^{n-r}
=-2|\scrx_r|m^{n-r-1}.
}
Then Claim~\ref{clnine} follows from the definition of $\phi_{n-r}$.
\claimend

\begin{clm}\label{clten}
For each $k\in \brk{n-r-1}$, we have 
\equ{eq5-15}
{
\sum_{A\in \scre\atop c(A)=k}
(-1)^{|A|} 
\left (|\scrs_A|-m^{c(A)}\right )
\ge -\phi_{k} m^{k}.
}
\end{clm}

\proofn
For each $A\in \scre$ with $c(A)=k$,
if $|A|$ is even,
\equ{eq5-13}
{
(-1)^{|A|} \left (|\scrs_A|-m^{c(A)}\right )
=|\scrs_A|-m^{c(A)}
\ge -m^{k}.
}
Thus Claim~\ref{clten} follows from Claim~\ref{clsix} 
and the definition of $\phi_k$.  \label{oc92}
\claimend

Let $\phi'_k$ be 
the number of subsets $A\subseteq E(G)$ 
such that $c(A)=k$, $\spnn{G}{A}$ is not a forest
and $|A|$ is even. 
Obviously, $\phi'_k\ge \phi_k$.
By the expression of (\ref{eq5-6}) 
and Claims~\ref{clfive},~\ref{clseven},~\ref{clnine} and~\ref{clten}, 
\eqn{eq5-14}
{
P_{DP}(G,\seth)-P(G,m)
&\ge &|\scrx_r|m^{n-r}
-2\phi_{n-r}|\scrx_r| m^{n-r-1}
-\sum_{k=1}^{n-r-1} \phi_{k}m^{k}
\nonumber \\
&\ge & m^{n-r}
-2\phi_{n-r} m^{n-r-1}
-\sum_{k=1}^{n-r-1} \phi_{k}m^{k} 
\nonumber \\
&\ge & m^{n-r}
-2\phi'_{n-r} m^{n-r-1}
-\sum_{k=1}^{n-r-1} \phi'_{k}m^{k}, 
}
where the second inequality holds 
when $m\ge 2\phi_{n-r}$.
As $\phi'_k$ is independent of the value of $m$,
by (\ref{eq5-14}), 
there must be a number $M_r\in \N$  
such that $P_{DP}(G,\seth)-P(G,m)> 0$ for 
all $m\ge M_r$.

Let $M=\max\{M_r: 3\le r\le n, r \mbox{ is odd}\}$.
Then, we conclude that for any $m\ge M$ and any $m$-fold cover 
$\seth=(L,H)$ of $G$,
if  $H\not\cong H_0(G,m)$, then
$P_{DP}(G,\seth)-P(G,m)> 0$ holds.

Hence Theorem~\ref{th1-2} holds. 
\proofend 

We end this section with an application
of Theorem~\ref{th1-2} to 
the generalized $\theta$-graphs.
For any $k$ numbers $a_1,a_2,\cdots,a_k\in \N$, 
where $k\ge 2$,
let $G=\Theta_{a_1,a_2,\cdots,a_k}$ 
denote the generalized $\theta$-graph obtained 
by connecting two distinct vertices with 
$k$ internally disjoint \label{cr29}
paths of lengths
$a_1,a_2,\cdots,a_k$ respectively.

Assume that $a_1\le a_2\le \cdots\le a_k$ 
and $a_1+a_2\ge 3$. 
Halberg, Kaul, Liu, Mudrock, Shin
and Thomason \cite{Halb}
showed that 
$P_{DP}(\Theta_{a_1,a_2,\cdots,a_k})
\approx P(\Theta_{a_1,a_2,\cdots,a_k})$
if $a_1+a_i$ is odd for each $i\in \brk{k}\setminus \{1\}$,
and 
$\dpl{\Theta_{a_1,a_2,\cdots,a_k}}$
otherwise. 

In the case that $a_1+a_i$ is odd for each 
$i\in \brk{k}\setminus \{1\}$,
$\Theta_{a_1,a_2,\cdots,a_k}$ belongs to the set $\gt_0$,
and thus $\Theta_{a_1,a_2,\cdots,a_k}\in \DPg$
by Theorem~\ref{th1-2},
implying that $\dpg{\Theta_{a_1,a_2,\cdots,a_k}}$.

\label{oc10}

\section{Proof of Theorem~\ref{th1-2-0}
\label{nsec5}}

For a chordal graph $G$,  
$P_{DP}(G,m)=P(G,m)$ for all $m\in \N$ (see \cite{kaul1}),
and thus, $\dpg{G}$ holds.
In the following, we first generalize this 
conclusion to some non-chordal graphs containing 
simplicial vertices.

\prop{pro2-4}
{
Let $u$ be a simplicial vertex of $G$.
For each $m\in \N$ with $m\ge d(u)+1$, if $P_{DP}(G-u,m)=P(G-u,m)$,
then   
$P_{DP}(G,m)=P(G,m)$.
}

\proof Assume that $P_{DP}(G-u,m)=P(G-u,m)$.
For any $m$-fold cover $\seth=(L,H)$ of $G$,  
\equ{eq2-4}
{
P_{DP}(G, \seth)\ge (m-d(u))P_{DP}(G-u, \seth'_u)
\ge 
(m-d(u))P_{DP}(G-u,m)=(m-d(u))P(G-u,m),
}
where $\seth'_u=(L',H')$ is the $m$-fold cover of $G-u$
with $L'(w)=L(w)$ for each $w\in V(G)\setminus \{u\}$ 
and $H'=H-L(u)$.
Thus, $P_{DP}(G,m)\ge (m-d(u))P(G-u,m)=P(G,m)$
by (\ref{eq2-1}).
On the other hand, 
$P_{DP}(G,m)\le P(G,m)$. 
Thus, the result follows.
\proofend

The first part of Theorem~\ref{th1-2-0}
follows from Proposition~\ref{pro2-4} directly.
In order to prove the second part of Theorem~\ref{th1-2-0},  
we need to introduce 
some preliminary results.

For any cover $\seth=(L,H)$ of $G$,
let $\seti(H)$ denote the set of independent sets $I$ in $H$
with $|I|=|V(G)|$. 
Thus, $P_{DP}(G,\seth)=|\seti(H)|$.
The {\it coloring number} of $G$, 
denoted by $col(G)$, is   \label{cr16}
the smallest integer $d$ 
for which there exists an ordering, 
$v_1, v_2,\cdots, v_n$ of the elements in $V(G)$, where $n=|V(G)|$, 
such that $|N_G(v_i)\cap \{v_1, v_2,\cdots, v_{i-1}\}|<d$
for each $i\in \brk{n}$.
Obviously, $\chi_{DP}(G)\le col(G)\le n$.
If $|L(v)|\ge col(G)$ for all $v\in V(G)$, 
then $\seti(H)\ne \emptyset$.

The following  fundamental property is important for 
the study of DP coloring.

\prop{pro2-0}
{
Let $\seth=(L,H)$ be a cover of $G$
with $|L(v)|\ge |V(G)|$ for each $v\in V(G)$.
Then, each independent set $A$ of $H$ is a subset of 
some set $I$ in $\seti(H)$.  
}

\proof If $A=\emptyset$, then the conclusion 
follows from the 
the fact that $|V(G)|\ge col(G)$.

Now assume that $A=\{(v_i,\pi_i): i\in \brk{k}\}$, where $k\ge 1$.
Clearly, $v_1,v_2,\cdots,v_k$ are pairwise distinct.
Let $\seth'=(L',H')$ be the cover of the subgraph 
$G'=G-\{v_i: i\in \brk{k}\}$,
where $L'(v)=L(v)\setminus N_H(A)$
for each $v\in V(G')$
and $H'$ is the subgraph of $H$ induced by 
$\bigcup_{v\in V(G')} L'(v)$.

Observe that $|L'(v)|\ge |L(v)|-k\ge |V(G')|$ 
for each $v\in V(G')$.
By the conclusion for $A=\emptyset$, 
there exists $I'\in \seti(H')$,
implying that $I=A\cup I'\in \seti(H)$.
\proofend \label{cr17}

By Proposition~\ref{pro2-0}, 
the following corollary is obtained.  

\corr{co2-0}
{
For any cover $\seth=(L,H)$ of $G$ 
with $|L(v)|\ge |V(G)|$ for each $v\in V(G)$, 
if $\seth'=(L,H')$ is a cover of $G$,
where $H'$ is obtained from $H$ by removing any edge 
in some set $E_H(L(v_1),L(v_2))$, where $v_1\ne v_2$,
then $P_{DP}(G,\seth')>P_{DP}(G,\seth)$.
}

For any $u\in V(G)$ and an $m$-fold cover 
$\seth=(L,H)$ of $G$,
let $\seth'_u=(L',H')$ be the cover of $G-u$,
where $H'=H-L(u)$ and  
$L'(v)=L(v)$ for each $v\in V(G)\setminus \{u\}$.
For any $I'\in \seti(H')$, let 
$$
\seti_H(I')=\{I'\cup \{(u,i)\}\in \seti(H): (u,i)\in L(u)\}.
$$
Obviously, for $m\ge d(u)$ and $I'\in \seti(H')$,
$|\seti_H(I')|\ge (m-d(u))$ holds,
implying that for $m>d(u)$, 
\eqn{eq2-20}
{
P_{DP}(G,\seth)
&=&|\seti(H)|=\sum_{I'\in \seti(H')}|\seti_H(I')|
\ge  \sum_{I'\in \seti(H')}(m-d(u)) \nonumber \\
&=&(m-d(u))|\seti(H')|=(m-d(u))P_{DP}(G-u,\seth'_u),
}
where $P_{DP}(G,\seth)>(m-d(u))P_{DP}(G-u,\seth'_u)$
if $|\seti_H(I')|>m-d(u)$ for some $I'\in \seti(H')$.

\prop{pro2-1}
{Let $\seth=(L,H)$ be an $m$-fold cover of $G$, 
where $m\ge |V(G)|$, and $u\in V(G)$. 
Then  $P_{DP}(G,\seth)\ge (m-d(u))P_{DP}(G-u,\seth'_u)$,
where the inequality is strict under each of the following
conditions:
\begin{enumerate}
\item $|E_H(L(u),L(v))|\le m-1$ for some $v\in N_G(u)$; or
\item $N_H((u,i))\setminus L(u)$ is not a clique of $H$
for some vertex $(u,i)\in L(u)$. 
\end{enumerate}
}

\proof By (\ref{eq2-20}), 
$P_{DP}(G,\seth)\ge (m-d(u))P_{DP}(G-u,\seth'_u)$ holds.
We need to prove that 
$P_{DP}(G,\seth)> (m-d(u))P_{DP}(G-u,\seth'_u)$ 
if either condition (i) or (ii) is satisfied.

Assume that condition (i) holds, 
i.e., $|E_H(L(u),L(v))|\le m-1$ for some $v\in N_G(u)$.
Then, there exists a $m$-fold cover
$\seth^*=(L, H^*)$  of $G$,
where $H^*$ is obtained from $H$ 
by adding a new edge joining some vertex in $L(u)$ to 
some vertex in $L(v)$.
By Corollary~\ref{co2-0}, 
\equ{}
{
P_{DP}(G,\seth)>P_{DP}(G,\seth^*)\ge (m-d(u))P_{DP}(G-u,\seth'_u).
}

Now assume that condition (ii) holds.
Without loss of generality, 
assume that 
$N_H((u,1))\setminus L(u)$ is not a clique of $H$.
Let $(v_1,i_1)$ and $(v_2,i_2)$ be 
non-adjacent vertices in $N_H((u,1))\setminus L(u)$.
Clearly, $v_1\ne v_2$.

As $\seth'_u=(L',H')$ is an $m$-fold cover of $G-u$
and $m\ge |V(G)|$, 
by Proposition~\ref{pro2-0},
there exists $I'\in \seti(H')$ such that 
$\{(v_1,i_1),(v_2,i_2)\}\subseteq I'$.

Note that $|I'\cap L(v)|=1$ for each $v\in N_G(u)$ and 
$\{(v_1,i_1),(v_2,i_2)\}\subseteq I'\cap N_H((u,1))$.
Assume that $N_G(u)=\{v_1,v_2,\cdots,v_r\}$, where $r=d(u)$, 
and $I'\cap L(v_j)=\{(v_j,\pi_j)\}$ for all $j\in [r]$.
Then
\eqn{eq40}
{
\left | L(u)\cap 
\bigcup_{j\in [r]}N_H((v_j,\pi_j)) \right |
&\le & \left |L(u)\cap\bigcup_{j\in \brk{2}}N_H((v_j,\pi_j)) \right |
+\left |L(u)\cap\bigcup_{3\le j\le r}N_H((v_j,\pi_j)) \right |
\nonumber \\
&\le & |\{(u,1)\}|+(r-2)=d(u)-1,
}
implying that 
$$
|\seti_H(I')|=m-\left | L(u)\cap 
\bigcup_{j\in [r]}N_H((v_j,\pi_j)) \right |
\ge m-d(u)+1.
$$

By (\ref{eq2-20}), 
$P_{DP}(G,\seth)>(m-d(u))P_{DP}(G-u,\seth'_u)$
holds.
The result is proven.
\proofend

We are now ready to prove Theorem~\ref{th1-2-0}
by applying (\ref{eq2-1}) and 
Propositions~\ref{pro2-4} and~\ref{pro2-1}.

\noindent {\it Proof of Theorem~\ref{th1-2-0}}:
If $\dpg{G-u}$, then
$\dpg{G}$ due to Proposition~\ref{pro2-4}.

Now assume that $G-u\in \DPg$.
Then, there exists $M\in \N$ such that
$P_{DP}(G-u,\seth')>P(G-u,m)$ for each integer $m\ge M$
and every $m$-fold cover $\seth'=(L',H')$ of $G-u$ with 
$H'\not\cong H_0(G-u,m)$.
 
Let $\seth=(L,H)$ be any $m$-fold cover of $G$ 
such that $H\not\cong H_0(G,m)$.
We may assume that 
$L(v)=\{(v,i): i\in \brk{m}\}$ for each $v\in V(G)$.
If $|E_H(L(v_1),L(v_2))|<m$ for some edge 
$v_1v_2\in E(G)$,
then, by Corollary~\ref{co2-0}, 
$P_{DP}(G, \seth)>P_{DP}(G, \seth^*)$ for $m\ge |V(G)|$, 
where $\seth^*$ is the $m$-fold cover $(L,H^*)$
obtained from $\seth$ by adding a new edge 
joining a vertex in $L(v_1)$ to a vertex in $L(v_2)$. 
Therefore, we can assume that 
$|E_H(L(v_1),L(v_2))|=m$ for each edge $v_1v_2\in E(G)$
and $H\not\cong H_0(G,m)$.

Consider the $m$-fold cover $\seth'_u=(L',H')$ of $G-u$.

\noindent {\bf Case 1}:
$H'\not \cong H_0(G-u,m)$.

By the assumption in the beginning of the proof,  
$P_{DP}(G-u,\seth'_u)>P(G-u,m)$ 
for each integer $m\ge M$.
By (\ref{eq2-1}) and Proposition~\ref{pro2-1},
for $m\ge \max\{M,|V(G)|\}$, 
\equ{pe2-3-1}
{
P_{DP}(G,\seth)\ge (m-d(u))P_{DP}(G-u,\seth')
>(m-d(u))P(G-u,m)
=P(G,m).
}

\noindent {\bf Case 2}:
$H'\cong H_0(G-u,m)$.

We can assume that $H'=H_0(G-u,m)$.
Since $H\not\cong H_0(G,m)$,
there must be some vertex 
$(u,i)\in L(u)$ that is \label{cr31}
adjacent to 
two vertices $(v_1,i_1)$ and $(v_2,i_2)$ 
with $v_1\ne v_2$ and $i_1\ne i_2$.
Since $H'=H_0(G-u,m)$ and $i_1\ne i_2$, 
$(v_1,i_1)$ and $(v_2,i_2)$ are not adjacent in $H$,
implying that $N_H((u,i))\setminus L(u)$ is not a clique of $H$.
By Proposition~\ref{pro2-1} again,
$P_{DP}(G,\seth)>(m-d(u))P(G,m)$ for $m\ge |V(G)|$.

Thus Theorem~\ref{th1-2-0} holds.
\proofend

By Theorem~\ref{th1-2-0},
we have the following consequence,
which generalizes the known conclusion that 
$\dpg{G}$ holds for every chordal graph $G$.

\begin{cor}\label{co4-5}
Let $G_1$ and $G_2$ be vertex-disjoint graphs 
and $k\in \N$, where $k\le \min\{\omega(G_i):i=1,2\}$.
Assume that $G_1$ is  chordal and $G\in {\mathscr G}(G_1\cup_k G_2)$.
If $\dpg{G_2}$, then 
$\dpg{G}$; also, if $G_2\in \DPg$, then  $G\in \DPg$.
\end{cor}

\proof 
As $G_1$ is chordal, there must be 
an ordering $v_1,v_2,\cdots,v_r$
of vertices in $V(G)\setminus V(G_2)$, where $r=|V(G_1)|-k$, 
such that 
$v_i$ is a simplicial vertex in 
$G-\{v_j: j\in \brk{i-1}\}$ for each $i\in \brk{r}$.
Then, the result follows  \label{cr32}
 from Theorem~\ref{th1-2-0}. 
\proofend

We wonder if Corollary~\ref{co4-5} holds 
without the condition that $G_1$ is chordal.

\prom{prob5}
{
For any vertex-disjoint graphs $G_1$ and $G_2$
and $k\in \N$, where $k\le \min\{\omega(G_i):i=1,2\}$,
is it true that 
if $\dpg{G_i}$ for $i=1,2$, then 
$\dpg{G}$ for every graph 
$G\in {\mathscr G}(G_1\cup_k G_2)$;
also, if $G_1,G_2\in \DPg$, 
then ${\mathscr G}(G_1\cup_k G_2)\subseteq \DPg$?
}

\section*{Acknowledgement}

The authors would like to thank the referees 
for their very helpful comments and suggestions.
The second author, Yan Yang, 
is supported by NSFC (No. 11971346).

\end{document}